\documentstyle[12pt]{article}
\def\i{{\infty}}
\def\L{{L^2(X;A)}}
\def\Hom{\mathop{\rm Hom}\nolimits}
\def\hx{{\chi_{\raisebox{-1pt}{\scriptsize $X$}}}}
\def\hj{{\chi_{\raisebox{-1pt}{\scriptsize $E_{n+1,j}$}}}}
\def\h{{\chi_{\raisebox{-1pt}{\scriptsize $E$}}}}
\def\hi{{\chi_{\raisebox{-1pt}{\scriptsize $E_i$}}}}
\date{12 January 1995}
\author{V.~M.~Manuilov}
\title{Lusin's $C$-property is not valid for functional Hilbert modules}

\begin{document}
\maketitle
\begin{abstract}
 We show that elements of Hilbert $A$-module obtained by completion
of the space of square-integrable functions on a space with measure $X$
taking values in a $C^*$-algebra $A$ cannot be viewed as $A$-valued functions
on $X$ defined almost everywhere
\end{abstract}

Let $A$ be a $C^*$-algebra.
Let $X$ be a space with finite measure $\mu$ and let $\lbrace E_i\rbrace$
be a countable
division of $X$ by non-intersecting measurable subsets. We call as usual a
function $f$ on $X$ {\it elementary}\/ if it is of the form
$$f(x)=\sum_i a_i\cdot \hi,$$
where $\h$ is the characteristic function of $E$.
Define an inner $A$-valued product on the space of elementary functions by
$$\langle f,g\rangle = \sum a_i^* b_i \,\mu(E_i),$$ where $\lbrace E_i
\rbrace$ is a
subdivision of divisions for $f$ and $g$. We call such function {\it
square-integrable}\/ if the series
\begin{equation}
\Vert f\Vert^2 =\langle f,f\rangle =\sum a_i^* a_i\,\mu(E_i)
\end{equation}
converges in norm in $A$. Completion of the space of square-integrable
elementary
functions with respect to the norm (1) we denote by $\L$.
Existance of the Walsh basis in $L^2(X)$ consisting of functions which are
equal to $\pm 1$ almost everywhere makes
this definition
equivalent to the standart one (cf \cite{ka1}, \cite{mf}),
where it is shown
that $\L$ is a Hilbert $A$-module isomorphic to the  Hilbert module
$l_2(A)$ \cite{pa}.

\bigskip
Unlike the case of the space of {\bf C}-valued square-integrable functions
it is impossible to consider elements of $\L$ as $A$-valued
functions on $X$ defined almost everywhere. Also there is nothing like the
classical Lusin's $C$-property. It is so even in the most
``convenient'' case of $A$ being a commutative $W^*$-algebra.

\bigskip
{\bf Example.} Let $X=Y=[0;1]$ with the standart measure $\mu$, and let
$A=L^{\i}(Y)$. Define a sequence of elementary functions $\lbrace f_n\rbrace
\in \L$
by induction. Put $f_1(x)=\hx(y)\in A$. On the $n$-th step divide $X$
into $2^{n-1}$ equal intervals $E_{n,i},\ i=1,\ldots,2^{n-1}$. Going to the
$(n{+}1)$-th step divide each of these intervals into two equal intervals
$E_{n+1,2i-1}$ and $E_{n+1,2i}$. Suppose that on the $n$-th step we have
$$f_n(x)=g_i(y)\in A,\quad {\rm when}\ x\in E_{n,i}$$
and define
$$ f_{n+1}(x)=g_i(y)+\hj(y),\quad
i={\rm E}({\textstyle\frac{j+1}{2}})\quad {\rm when}\ x\in E_{n+1,j}.$$
Obviously we have
$$\Vert f_n -f_{n+1}\Vert_{\L}=\sum_{j=1}^{2^n}
\hj(y)\cdot\mu(E_{n+1,j})=\frac{1}{2^n}\sum_{j=1}^{2^n}\hj(y)=
\frac{1}{2^n},$$ hence the sequence
$\lbrace f_n\rbrace$ is fundamental in $\L$. But for almost all $x\in X$
$$\mu(\lbrace y\in Y:\ f_n(x)\geq n\rbrace)=\frac{1}{2^{n-1}},$$
so almost everywhere on $X$ $f(x)\in A$ fails.

\bigskip
{\bf Remark.} Nevertheless  for any $C^*$-algebra $A$ and
for any positive functional $\phi\in A^*$ we can prolonge its action from
$A$ to elementary functions by
$$\bar\phi\Bigl(\sum a_i\hi(x)\Bigr)=\sum\phi(a_i)\hi(x)$$
and then to
$\L$ and even to its dual module $\L'=\Hom_A(\L;A)$,
$\bar\phi : \L\longrightarrow L^2(X)$
and the estimate
$$
\Vert\bar\phi(f)\Vert^2_{L^2(X)}=\sum\vert\phi(a_i)\vert^2\,\mu(E_i)
\leq\sum\phi(a^*_ia_i)\,\mu(E_i) =
\phi(\langle f,f\rangle)\leq\Vert f\Vert_{\L'}^2
$$
for $f = \sum a^*_ia_i\hi(x)$
shows that $\bar\phi(f)$
is square-integrable
on $X$ for any $f\in \L'$ and for any positive $\phi\in A^*$.

\bigskip
{\bf Acknoledgement.\/}
This work was partially
supported by the Russian Foundation for Fundamental Research (grant
\mbox{N 94-01-00108-a)} and the International Science Foundation
(grant N MGM000).
I am indebted to M.~Frank, A.~A.~Irmatov, A.~S.~Mishchenko
and E.~V.~Troitsky for helpful discussions.

\vspace{2cm}
V.M.Manuilov \\*
Moscow State  University \\*
Russia - 119899 Moscow \\*
E-mail:manuilov@math.math.msu.su

\end{document}